\documentclass[11pt]{amsart}
\usepackage{amsmath,amsthm,amssymb, amscd, amsfonts}
\usepackage[all]{xy}
\usepackage{leftidx}
\usepackage[latin1]{inputenc}
\usepackage{enumerate}
\usepackage[dvipsnames]{xcolor}

\theoremstyle{plain}

\newtheorem{theorem}{Theorem}[section]

\newtheorem{proposition}[theorem]{Proposition}

\theoremstyle{definition}

\theoremstyle{remark}

\numberwithin{equation}{section}\theoremstyle{plain}

\newcommand{\I}{\mathcal{I}}

\renewcommand{\1}{\textbf{1}}

\newcommand{\C}{{\mathcal C}}
\newcommand{\D}{{\mathcal D}}
\newcommand{\F}{{\mathcal F}}

\newcommand{\Z}{{\mathcal Z}}

\newcommand\Irr{\operatorname{Irr}}
\newcommand\FPdim{\operatorname{FPdim}}

\newcommand\vect{\operatorname{Vec}}

\begin{document}
\title[Rank $4$ fusion categories]{On Grothedieck rings of rank $4$ self-dual fusion categories}
\author[J. Dong]{Jingcheng Dong$^{1,2}$}
\email{jcdong@nuist.edu.cn}
\address{1. College of Mathematics and Statistics, Nanjing University of Information Science and Technology, Nanjing 210044, China}
\address{2. Center for Applied Mathematics of Jiangsu Province, Nanjing University of Information Science and Technology, Nanjing 210044, China}



\date{\today}

\begin{abstract}
Let $\C$ be a self-dual fusion category of rank $4$ which has a nontrivial proper fusion subcategory. We identify three new families of Grothendieck rings for $\C$: one of them is completely determined, the other two are parameterized by several non-negative integers.
\end{abstract}

\maketitle

\section{Introduction}\label{sec1}
The study of fusion categories, which are semisimple rigid tensor categories with finitely many simple objects and finite-dimensional morphism spaces, has been a central topic in the field of mathematical physics and representation theory. These categories generalize the notion of group representations and have found applications in various areas, including conformal field theory, topological quantum field theory, and subfactor theory.

The classification of fusion categories is a challenging problem, and significant progress has been made primarily for categories of low rank. The rank of a fusion category, defined as the number of isomorphism classes of its simple objects, serves as a measure of its complexity. The classification of fusion categories of rank $1$ is straightforward, as there is only one such category-the category of vector spaces over a field.

The classification of fusion categories of rank $2$ was achieved by Ostrik in his seminal work \cite{2013ostrikpivotal}. Ostrik showed that there are exactly four fusion categories of rank 2, characterized by their Grothendieck rings, which are isomorphic to one of two based rings, $K_0$ and $K_1$. The Grothendieck ring $K_0$ corresponds to the category of representations of the cyclic group $\mathbb{Z}_2$, while $K_1$ is associated with the Yang-Lee fusion rules and is well-known in conformal field theory. Ostrik's work laid the foundation for understanding the structure of low-rank fusion categories and highlighted the importance of the Grothendieck ring as a fundamental invariant.

Building on Ostrik's results, subsequent research focused on fusion categories of rank $3$. In \cite{ostrik2003fusion}, Ostrik classified fusion categories of rank $3$ that admit a pivotal structure, a mild assumption that is conjectured to hold for all fusion categories. This classification revealed a rich structure involving pointed categories, categories associated with quantum groups, and Ising categories. The work of Ostrik and others in this area demonstrated the intricate interplay between fusion categories and other areas of mathematics, such as modular tensor categories and subfactor theory.

The classification of fusion categories of rank $4$ remains an open problem, with significant progress made in understanding specific subclasses. For instance, the work of Larson \cite{Larson2014184} focuses on pseudo-unitary non-self-dual fusion categories of rank $4$, providing a nearly complete classification of those based rings that admit pseudo-unitary categorification. Larson's study highlights the complexity of fusion categories of rank $4$ and the challenges involved in their classification. In \cite{Dong2017Non}, the authors complete the classification of the Grothendieck ring $K(C)$ of self-dual spherical fusion categories of rank $4$ with non-trivial grading, proving that $K(C) \cong \text{Fib} \otimes \mathbb{Z}[\mathbb{Z}_2]$, where Fib is the Fibonacci fusion ring and $\mathbb{Z}[\mathbb{Z}_2]$ is the group ring on $\mathbb{Z}_2$. They further show that if the category is braided, it is equivalent to $\text{Fib} \boxtimes \text{Vec}_{\mathbb{Z}_2}^{\omega}$ as fusion categories. 

Recently, Michell, Izumi and Penneys completes the classification of rank $4$ unitary fusion categories (UFC) that are $\mathbb{Z}_2$-quadratic and possess a dual pair of simple objects \cite{MIP2024}. The main tools for this classification are skein theory, a generalisation of Ostrik's results on formal codegrees. Rank $4$ UFCs also appear in the work of Liu, Palcoux and Wu \cite{LIU2021107905}, as an example to  illustrate the application of the Schur product criterion, which provides a necessary condition for unitary categorification.

We delve into the structure of self-dual fusion categories of rank $4$, considering cases where a nontrivial proper fusion subcategory has rank $2$ or $3$. We identify three new families of Grothendieck rings for self-dual fusion categories of rank 4, each parameterized by several non-negative integers.


The main aim of this paper is to construct some preliminary results and thus allow other people to do more detailed analysis. Throughout this paper, we shall work over an algebraically closed field $k$ of characteristic $0$. We refer to \cite{egno2015} for the main notions about fusion categories.

\section{Preliminaries}\label{prels}
Consider a fusion category $\C$. Denote by $\Irr(\C)$ the collection of equivalence classes of simple objects within $\C$. The set $\Irr(\C)$ forms a $\mathbb{Z}^+$ basis for the Grothendieck ring $K_0(\C)$. 

For any objects $X,Y$ belonging to $\Irr(\C)$, we have 
$$X\otimes Y=\bigoplus_{Z\in \Irr(\C)}N_{XY}^ZZ,$$
where $N_{XY}^Z$, called the fusion coefficient, is the multiplicity of $Z$ in the decomposition of $X\otimes Y$. The fusion coefficients obey the Frobenius Reciprocity as follows. For all $X,Y,Z\in \Irr(\C)$,
$$N_{XY}^Z = N_{ZY}^{X*}=N_{ZX}^{Y*}=N_{Y*X}^Z = N_{YZ}^{X*}=N_{X*Z}^Y.$$

The Frobenius-Perron dimension $\FPdim(X)$ of $X$ is the largest positive eigenvalue of the fusion matrix $N_X$ with entries $(N_{XY}^Z)$.  Such an eigenvalue exists by the Frobenius-Perron theorem. By \cite[Theorem 8.6]{etingof2005fusion}, the Frobenius-Perron dimension extends to a ring homomorphism $\FPdim : K_0(\C)\to \mathbb{R}$.

The formal codegrees of the fusion category $\C$ are the eigenvalues of the matrix 
$$M=\sum_{X\in \Irr(\C)}N_ZN_{Z^*}.$$

Let $\C$ be a fusion category. The Drinfeld  center $\Z(\C)$ of $\C$  is defined as the category whose objects are pairs $(X, c_{-,X} )$, where $X$ is an object of $\C$ and $c_{-,X}$ is a natural family of isomorphisms $c_{V,X} : V\otimes X\to X\otimes V$, $V\in \C$, satisfying certain compatibility conditions, see \cite[Definition XIII.4.1]{kassel1995quantum}.  It is shown in \cite[Theorem 2.15, Proposition 8.12]{etingof2005fusion} that $\Z(\C)$ is a braided fusion category and $\FPdim(\Z(\C))=\FPdim(\C)^2$. In addition, $\Z(\C)$ is non-degenerate by \cite[Corollary 3.9]{drinfeld2010braided}.

There is an obvious forgetful tensor functor $F:\Z(\C)\to \C$. By \cite[Proposition 3.39]{etingof2004finite}, the forgetful functor $F:\Z(\C)\to \C$ is surjective. Since $\Z(\C)$ is semisimple the forgetful functor $F$ has a right adjoint $I: \C\to \Z(\C)$. By \cite[Lemma 3.2]{etingof2011weakly}, $I(\1)\in \Z(\C)$ has a natural structure of commutative algebra. By \cite[Proposition 5.4]{etingof2005fusion}, $F(I(X))=\oplus_{V\in \Irr(\C)}V\otimes X\otimes V^*$.

\section{The Grothendieck rings}\label{sec3}
Let $\C$ be a self-dual fusion category of rank $4$. Assume that $\C$ contains a nontrivial proper fusion subcategory $\D$. If $\D$ has rank $3$ then $\C$ is a near-group fusion category $\C(\mathbb{Z}_3,k)$ \cite[Theorem 6.2]{DongChenWang2022}, where $k=0,2,3$ or $6$. The Grothendieck ring of such a fusion category is well known. In the rest of this paper, we only consider the cases when $\D$ has rank $2$. By the work of Ostrik \cite{ostrik2003module}, $\D$ is either pointed or a Fibonacci fusion category. In addition, we assume that $\C$ is not pointed since self-dual pointed fusion category of rank $4$ is the category $\vect_{\mathbb{Z}_2\times\mathbb{Z}_2}^{\omega}$  of $G$-graded vector spaces with associativity constraint given by a $3$-cocycle $\omega\in H^3(\mathbb{Z}_2\times\mathbb{Z}_2,k^{\times})$.

\subsection{$\D$ is pointed}\label{sec3.1}
\begin{proposition}\label{prop31}
Let $\C$ be a self-dual fusion category of rank $4$. Assume that $\Irr(\C)=\{\1,g,X,Y\}$ and $\{\1,g\}$ generates $\C_{pt}$. If the action $g\otimes-$ on $\Irr(\C)/G(\C)$ is fixed-point free then $\C$ has one of the following fusion rules:
\begin{equation}\label{eq01}
\begin{split}
&g\otimes g=\1,g\otimes X=X, g\otimes Y=Y,X\otimes X=\1\oplus mX\oplus mY,\\
&Y\otimes X=g\oplus mX\oplus mY, Y\otimes Y=\1\oplus mX\oplus mY, m=1, 2.
\end{split}
\end{equation}

\begin{equation}\label{eq02}
\begin{split}
&g\otimes g=\1,g\otimes X=Y, g\otimes Y=X,\\
X\otimes X&=\1\oplus X,Y\otimes X=g\oplus Y, Y\otimes Y=\1\oplus X.
\end{split}
\end{equation}
Moreover, (\ref{eq02})  is realized by $\rm Fib\otimes \vect_{\mathbb{Z}_2}^{\omega}$ for some $3$-cocycle $\omega\in H^3(\mathbb{Z}_2,k^{\times})$.
\end{proposition}
\begin{proof}
Since the action of $g$ is fixed-point free, we get $g\otimes X=Y, g\otimes Y=X$. From $N_{XX}^g=N_{gX}^X=N_{Y}^X=0$  and $N_{YY}^g=N_{gY}^Y=N_{X}^Y=0$ we know that the decomposition of $X\otimes X$ and $Y\otimes Y$ can not contain $g$ as a summand. Hence we may assume that $X\otimes X=\1\oplus mX\oplus nY$, where $m,n\in \mathbb{Z}^+$. From $N_{XX}^Y=N_{YX}^X=n$ and $N_{gX}^Y=N_{YX}^g=1$, we may write $Y\otimes X=g\oplus nX\oplus sY$, where $s\in \mathbb{Z}^+$. From $N_{YX}^Y=N_{XY}^Y=N_{YY}^X=s$, we may write $Y\otimes Y=\1\oplus sX\oplus tY$, where $t\in \mathbb{Z}^+$.

Applying  $\FPdim$ on both sides of the decomposition of $X\otimes X$, $Y\otimes Y$ and $Y\otimes X$, and using the fact $\FPdim(X)=\FPdim(Y)$, we get
$$m+n=s+t=n+s.$$
We thus have $m=s$ and $n=t$. Without loss of generality, we may assume that $m\geq n$.

It follows that the fusion matrices $M_g$, $M_X$ and $M_Y$ are
$$	M_{g}=\left(
\begin{matrix}
		0 & 1 & 0 & 0 \\
		1 & 0 & 0 & 0 \\
		0 & 0 & 0 & 1 \\
		0 & 0 & 1 & 0
	\end{matrix}
\right),
	M_{X}=\left(
\begin{matrix}
		0 & 0 & 1 & 0 \\
		0 & 0 & 0 & 1 \\
		1 & 0 & m & n \\
		0 & 1 & n & m
	\end{matrix}
\right),
	M_{Y}=\left(
\begin{matrix}
		0 & 0 & 0 & 1 \\
		0 & 0 & 1 & 0 \\
		0 & 1 & n & m \\
		1 & 0 & m & n
	\end{matrix}
\right).$$
Then $\FPdim(X)=\FPdim(Y)=(m+n+\sqrt{4+(m+n)^2})/2$.

Let $A=I_4+M_g^2+M_X^2+M_Y^2$. The the formal codegrees of $\C$ are
\begin{equation}\label{eq2}
\begin{split}
&f_1=4+(m+n)^2+(m+n)\sqrt{4+(m+n)^2},\\&f_2=4+(m+n)^2-(m+n)\sqrt{4+(m+n)^2},\\
&f_3=4+(m-n)^2+(m-n)\sqrt{4+(m-n)^2},\\&f_4=4+(m-n)^2-(m-n)\sqrt{4+(m-n)^2}.
\end{split}
\end{equation}
If $\sqrt{4+(m+n)^2}$ is a rational number then $m=n=0$. In this case, $\C$ is pointed, a contradiction. In the rest of our proof, we assume that $\sqrt{4+(m+n)^2}$ is an irrational number.

Since $K(\C)$ is a commutative ring, it has four $1$-dimensional representations. Hence the object $I(\1)$ is a direct sum of four simple objects which all have multiplicity $1$, see \cite[Theorem 2.13]{2013ostrikpivotal}. Therefore, we can assume $\I (\1) = \1 \oplus A \oplus B\oplus C$. Furthermore, we may assume, also by \cite[Theorem 2.13]{2013ostrikpivotal},
$$\FPdim(A)=\frac{f_1}{f_2},\FPdim(B)=\frac{f_1}{f_3},\FPdim(C)=\frac{f_1}{f_4}.$$
By \cite[Proposition 5.4]{etingof2005fusion}, we have
\begin{equation}\label{eq3}
\begin{split}
\F(\I(\1)) &= F(\1) \oplus \F(A) \oplus \F(B) \oplus \F(C)\\
 &= \mathop \oplus
\limits_{T \in \mbox{Irr}(C)} T \otimes T^\ast\\
& = 4 \cdot \1 \oplus 2mX \oplus 2nY.
\end{split}
\end{equation}
We assume that $\F(A)=\1 \oplus a_1X \oplus b_1Y$, $\F(B)=\1 \oplus a_2X \oplus b_2Y$, $\F(C)=\1 \oplus a_3X \oplus b_3Y$. From $\FPdim(\F(B))=1+(a_2+b_2)\FPdim(X)$, we get
\begin{equation}\label{eq4}
\begin{split}
\frac{f_1}{f_3}=1+\frac{a_2+b_2}{2}(m+n+\sqrt{4+(m+n)^2}).
\end{split}
\end{equation}

Set $\alpha=\sqrt{4+(m+n)^2}$, $\beta=m+n$, $\gamma=m-n$ and $\delta=\sqrt{4+(m-n)^2}$. Then the equation (\ref{eq4}) is reduced to the following equality:
\begin{equation}\label{eq5}
\begin{split}
(4+\gamma^2)[\beta^2-2\beta(a_2+b_2)]=[(4+\gamma^2)(2(a_2+b_2)-\beta)]\alpha+\gamma(4+\beta^2)\delta+\beta\gamma\alpha\delta.
\end{split}
\end{equation}
Observe that the left hand side of equality (\ref{eq5}) is a rational number. If $\delta$ is an irrational number and $\alpha\neq \delta$ then $\alpha,\delta$ and $\alpha\delta$ are different irrational numbers. Hence the right hand side of equality (\ref{eq5}) must be $0$. It follows that the coefficients of $\alpha,\delta$ and $\alpha\delta$ in  equality (\ref{eq5}) are $0$, which implies that $m=n$. But $\delta=2$ in this case which is a contradiction. Thus we have two possibilities:

(1)\, $\delta$ is a rational number.

(2)\, $\delta$ is an irrational number and $\alpha=\delta$.

If (1) holds true then $m=n$ and hence $m=1$ or $2$ by \cite[Theorem 2.4]{MIP2024}. This gets fusion rules (\ref{eq01}).

If (2) holds true then $n=0$ which means $X\otimes X=1\oplus mX$. Our assumption requires that $m\neq 0$ otherwise $\C$ is pointed. Hence $X$ generates a non-pointed fusion subcategory of rank $2$. By the classification of fusion category of rank $2$ \cite{ostrik2003module}, we get $m=1$. Thus we get fusion rules (\ref{eq02}).
\end{proof}

\begin{proposition}\label{prop32}
Let $\C$ be a self-dual fusion category of rank $4$. Assume that $\Irr(\C)=\{\1,g,X,Y\}$ and $\{\1,g\}$ generates $\C_{pt}$. If the action $g\otimes-$ on $\Irr(\C)/G(\C)$ has a fixed point then the fusion rules of $\C$ are determined by non-negative integers $m,n,s$ and $t$ as follows.
\begin{equation}\label{eq6}
\begin{split}
&g\otimes g=\1,g\otimes X=X, g\otimes Y=Y,Y\otimes X=nX\oplus sY,\\
X&\otimes X=\1\oplus g\oplus mX\oplus nY, Y\otimes Y=\1\oplus g\oplus sX\oplus tY.
\end{split}
\end{equation}
Moreover, $m,n,s$ and $t$ satisfy the following relation:
$$n^2+s^2=2+ms+nt.$$
\end{proposition}
\begin{proof}
Since the action $g\otimes-$ on $\Irr(\C)/G(\C)$ has a fixed point and $X,Y$ are the only simple objects except $\1$ and $g$, we get $g\otimes X=X$ and $g\otimes Y=Y$.

From $N_{gX}^X=N_{XX}^g=1$, we may write $X\otimes X=\1\oplus g\oplus mX\oplus nY$, where $m,n\in \mathbb{Z}^+$. From $N_{XX}^Y=N_{YX}^X=n$, we may write $Y\otimes X=nX\oplus sY$, where $s\in \mathbb{Z}^+$. From $N_{YX}^Y=N_{YY}^X=s$, we may write $Y\otimes Y=\1\oplus g\oplus sX\oplus tY$, where $s\in \mathbb{Z}^+$. We hence get the fusion rules as described above.

The fusion matrices $M_g$, $M_X$ and $M_Y$ are
$$	M_{g}=\left(
\begin{matrix}
		0 & 1 & 0 & 0 \\
		1 & 0 & 0 & 0 \\
		0 & 0 & 1 & 0 \\
		0 & 0 & 0 & 1
	\end{matrix}
\right),
	M_{X}=\left(
\begin{matrix}
		0 & 0 & 1 & 0 \\
		0 & 0 &1 & 0 \\
		1 & 1 & m & n \\
		0 & 0 & n &s
	\end{matrix}
\right),
	M_{Y}=\left(
\begin{matrix}
		0 & 0 & 0 & 1 \\
		0 & 0 & 0 & 1 \\
		0 & 0 & n & s \\
		1 & 1 & s & t
	\end{matrix}
\right).$$

The commutativity of $K(\C)$ implies these three matrices pairwise commute, which gives us the relation $n^2+s^2=2+ms+nt$.

\end{proof}

\subsection{$\D$ is a Fibonacci fusion category}\label{sec3.2}
\begin{proposition}\label{prop33}
Let $\C$ be a self-dual fusion category of rank $4$. Assume that $\Irr(\C)=\{\1,X,Y,Z\}$ and $\{\1,Z\}$ generates a Fibonacci fusion category. Then the fusion rules of $\C$ are determined by non-negative integers $a$ and $b$ as follows.
\begin{equation}\label{eq7}
\begin{split}
&X\otimes X=\1\oplus aX\oplus bY,X\otimes Y=bX\oplus (a+b)Y\oplus Z,X\otimes Z=Y,\\
Y\otimes &Z=X\oplus Y,Y\otimes Y=\1\oplus (a+b)X\oplus (a+2b)Y\oplus Z, Z\otimes Z=\1\oplus Z.
\end{split}
\end{equation}
\end{proposition}
\begin{proof}
Assume that $X\otimes X=\1\oplus aX\oplus bY\oplus cZ$, where $a,b,c\in \mathbb{Z}^+$. From $N_{XX}^Y=N_{YX}^X=b$, we may write $Y\otimes X=bX\oplus dY\oplus eZ$, where $d,e\in \mathbb{Z}^+$. Because $Z\otimes Z=\1\oplus Z$, the decomposition of $Z\otimes X$ and $Z\otimes Y$ can not contain $Z$ as summands. From $N_{YX}^Y=N_{YY}^X=d$, we may write $Y\otimes Y=\1\oplus dX\oplus gY\oplus fZ$, where $g,f\in \mathbb{Z}^+$. From $N_{YY}^Z=N_{ZY}^Y=f$ and $N_{YX}^Z=N_{ZX}^Y=N_{YZ}^X=e$, we get $Z\otimes Y=eX\oplus fY$. From $N_{XX}^Z=N_{ZX}^X=c$ and $N_{ZX}^Y=e$, we get $Z\otimes X=cX\oplus eY$.

The fusion matrices $M_X$, $M_Y$ and $M_Z$ are
$$	M_{X}=\left(
\begin{matrix}
		0 & 1 & 0 & 0 \\
		1 & a & b & c \\
		0 & b & d & e \\
		0 & c & e & 0
	\end{matrix}
\right),
	M_{Y}=\left(
\begin{matrix}
		0 & 0 & 1 & 0 \\
		0 & b & d & e \\
		1 & d & g & f \\
		0 & e & f & 0
	\end{matrix}
\right),
	M_{Z}=\left(
\begin{matrix}
		0 & 0 & 0 & 1 \\
		0 & c & e & 0 \\
		0 & e & f & 0 \\
		1 & 0 & 0 & 1
	\end{matrix}
\right).$$
The characteristic polynomial of $M_Z$ is
$$(x-\frac{1+\sqrt{5}}{2})(x-\frac{1-\sqrt{5}}{2})(x-\frac{c+f-\sqrt{(c-f)^2+4e^2}}{2})(x-\frac{c+f+\sqrt{(c-f)^2+4e^2}}{2})$$
Since $\FPdim(Z)=\frac{1+\sqrt{5}}{2}$ is the largest eigenvalue of $M_Z$, we have
\begin{equation}\label{eq8}
\begin{split}
c+f+\sqrt{(c-f)^2+4e^2}\leq 1+\sqrt{5}.
\end{split}
\end{equation}
The commutativity of $K(\C)$ implies that $M_XM_Y=M_YM_X$, $M_XM_Z=M_ZM_X$ and $M_YM_Z=M_ZM_Y$. Hence we get

\begin{equation}\label{eq9}
	\left\{
	\begin{aligned}
		& 1+ad+bg+cf=b^2+d^2+e^2 , & (1)\\
		& ae+bf=bc+de , & (2)\\
& be+df=cd+eg , & (3)\\
& c^2+e^2=1+c , & (4)\\
& e=e(c+f) , & (5)\\
& 1+f=e^2+f^2. & (6)
	\end{aligned}
	\right.
\end{equation}
By (4) and (6), we get $c^2+e^2\neq 0$ and $e^2+f^2\neq 0$. Hence inequality (\ref{eq8}) shows that we have three possibilities:

\begin{equation}\label{eq10}
	\left\{
	\begin{aligned}
		& c=1\\
		& e=0\\
& f=1
	\end{aligned}
	\right.;
\left\{
	\begin{aligned}
		& c=0\\
		& e=1\\
& f=1
	\end{aligned}
	\right.;
\left\{
	\begin{aligned}
		& c=1\\
		& e=1\\
& f=0
	\end{aligned}
	\right..
\end{equation}

If $c=1,e=0,f=1$ then $Z\otimes X=X$. This is impossible since $Z$ is not invertible.

If $c=0,e=1,f=1$ then we have the following from (1)-(6)
\begin{equation}\label{eq11}
	\left\{
	\begin{aligned}
		& ad+bg=b^2+d^2\\
		& a+b=d\\
        & b+d=g
	\end{aligned}
	\right.
\end{equation}
Further, $g=a+2b$. Hence the fusion rules of $\C$ are determined by $a$ and $b$ as follows.
\begin{equation*}\label{eq12}
\begin{split}
X\otimes X&=\1\oplus aX\oplus bY,X\otimes Y=bX\oplus (a+b)Y\oplus Z,X\otimes Z=X\oplus Y,\\
Y\otimes &Z=X\oplus Y,Y\otimes Y=\1\oplus (a+b)X\oplus (a+2b)Y\oplus Z, Z\otimes Z=\1\oplus Z.
\end{split}
\end{equation*}

If $c=1,e=1,f=0$ then similar arguments shows the fusion rules of $\C$ are determined by $d$ and $g$ as follows.
\begin{equation*}\label{eq13}
\begin{split}
X&\otimes X=\1\oplus (2d+g)X\oplus (d+g)bY,X\otimes Y=(d+g)X\oplus dY\oplus Z,Y\otimes Z=X,\\
&X\otimes Z=X\oplus Y,Y\otimes Y=\1\oplus (a+b)X\oplus (a+2b)Y\oplus Z, Z\otimes Z=\1\oplus Z.
\end{split}
\end{equation*}
It is clear this Grothendieck ring is isomorphic to the previous one.
\end{proof}


\begin{thebibliography}{10}

\bibitem{DongChenWang2022}
Jingcheng Dong, Gang Chen, and Zhihua Wang.
\newblock Fusion categories containing a fusion subcategory with maximal rank.
\newblock {\em J. Algebra}, 604(1):107--127, 2022.

\bibitem{Dong2017Non}
Jingcheng Dong, Liangyun Zhang, and Li~Dai.
\newblock Non-trivially graded self-dual fusion categories of rank $4$.
\newblock {\em Acta Math.Sin.}, 34(2):275--287, 2017.

\bibitem{drinfeld2010braided}
Vladimir Drinfeld, Shlomo Gelaki, Dmitri Nikshych, and Victor Ostrik.
\newblock On braided fusion categories {I}.
\newblock {\em Selecta Math., New Ser.}, 16(1):1--119, 2010.

\bibitem{MIP2024}
Cain Edie-Michell, Masaki Izumi, and David Penneys.
\newblock Classification of $\mathbb{Z}/2\mathbb{Z}$-quadratic unitary fusion
  categories.
\newblock {\em Quantum Topol}, DOI: 10.4171/QT/201, 2024.

\bibitem{egno2015}
P.~Etingof, S.~Gelaki, D.~Nikshych, and V.~Ostrik.
\newblock {\em Tensor Categories}.
\newblock Mathematical surveys and monographs, vol. 205. Amer. Math. Soc.,
  2015.

\bibitem{etingof2011weakly}
Pavel Etingof, Dmitri Nikshych, and Victor Ostrik.
\newblock Weakly group-theoretical and solvable fusion categories.
\newblock {\em Adv. Math.}, 226(1):176--205, 2011.

\bibitem{etingof2005fusion}
Pavel Etingof, Dmitri Nikshych, and Viktor Ostrik.
\newblock On fusion categories.
\newblock {\em Ann. Math.}, 162(2):581--642, 2005.

\bibitem{etingof2004finite}
Pavel Etingof and Viktor Ostrik.
\newblock Finite tensor categories.
\newblock {\em Mosc. Math. J}, 4(3):627--654, 2004.

\bibitem{kassel1995quantum}
Christian Kassel.
\newblock Quantum groups, {GTM} 155, 1995.

\bibitem{Larson2014184}
Hannah~K. Larson.
\newblock Pseudo-unitary non-self-dual fusion categories of rank 4.
\newblock {\em J. Algebra}, 415:184 -- 213, 2014.

\bibitem{LIU2021107905}
Zhengwei Liu, Sebastien Palcoux, and Jinsong Wu.
\newblock Fusion bialgebras and {F}ourier analysis: Analytic obstructions for
  unitary categorification.
\newblock {\em Adv. Math.}, 390:107905, 2021.

\bibitem{ostrik2003module}
Victor Ostrik.
\newblock Module categories, weak {H}opf algebras and modular invariants.
\newblock {\em Transform. Groups}, 8(2):177--206, 2003.

\bibitem{2013ostrikpivotal}
Victor Ostrik.
\newblock Pivotal fusion categories of rank 3.
\newblock {\em Mosc. Math. J.}, 15(2):373--396, 2015.

\bibitem{ostrik2003fusion}
Viktor Ostrik.
\newblock Fusion categories of rank 2.
\newblock {\em Math. Res. Lett.}, 10(2):177--183, 2003.

\end{thebibliography}
\end{document}